\documentclass[12pt,a4paper]{article}

\usepackage{amsmath,amssymb,amsthm}

\newtheorem{theorem}{Theorem}
\newtheorem{lemma}[theorem]{Lemma}

\newtheorem{corollary}[theorem]{Corollary}

\newtheorem*{theorema}{Theorem A}
\newtheorem*{theoremb}{Theorem B}
\newtheorem*{theoremc}{Theorem C}

\newtheorem*{hypothesis}{Hypothesis}

\newcommand{\spaceP}{\mathcal{P}}

\begin{document}
\pagestyle{plain}
\pagenumbering{arabic}
\title{Nilpotent Singer groups}
\author{Nick Gill}
\maketitle

\begin{abstract}
Let $N$ be a nilpotent group normal in a group $G$. Suppose that $G$ acts transitively upon the points of a finite non-Desarguesian projective plane $\spaceP$. We prove that, if $\spaceP$ has square order, then $N$ must act semi-regularly on $\spaceP$. 

In addition we prove that if a finite non-Desarguesian projective plane $\spaceP$ admits more than one nilpotent group which is regular on the points of $\spaceP$ then  $\spaceP$ has non-square order and the automorphism group of $\spaceP$ has odd order.

\end{abstract}

\section{Introduction}

A {\it Singer group} $S$ of a projective plane $\spaceP$ of order $x$ is a collineation group of $\spaceP$ which acts sharply transitively on the points of $\spaceP$. The existence of such a Singer group is equivalent to a $(v,k,1)$ difference set in $S$ where $v=x^2+x+1$ and $k=x+1$.

Ho \cite[theorem 1]{ho} has proved the following theorem concerning {\it abelian} Singer groups:

\begin{theoremc}
A finite projective plane which admits more than one abelian Singer group is Desarguesian.
\end{theoremc}

We will present an alternative proof of this theorem (our proof, unlike Ho's, will be dependent on the Classification of Finite Simple Groups) and then will present work aimed at extending the result to {\it nilpotent} Singer groups. In particular we prove the following:

\begin{theoremb}
Suppose that a non-Desarguesian finite projective plane $\spaceP$ of order $x$ admits more than one nilpotent Singer group. Then the automorphism group of $\spaceP$ has odd order and $x$ is not a square.
\end{theoremb}
 
In the course of proving Theorem B we will need to prove the following:

\begin{theorema}\label{theorem:nilpotentnormal}
Let $F$ be a nilpotent group which is normal in a transitive automorphism group $G$ of $\spaceP$, a projective plane of order $x=u^2$. Then $F$ acts semi-regularly on $\spaceP$. 
\end{theorema}

For the remainder of the paper we operate under the following hypothesis. The conditions included represent, by \cite{wagner} and \cite[4.1.7]{dembov}, the conditions under which a group may act transitively on the points of a non-Desarguesian projective plane.

\begin{hypothesis}
Let $\spaceP$ be a non-Desarguesian projective plane of order $x>4$. Let $G$ be an automorphism group of $\spaceP$ which acts transitively upon the points of $\spaceP$. If $G$ contains any involutions then $x=u^2$, $u>2$, and each involution fixes $u^2+u+1$ points. If $h\in G$ then $h$ fixes at most $x+\sqrt{x}+1$ points.
\end{hypothesis}

By \cite[Theorem A]{gill}, we know that the Fitting group and the generalized Fitting group of $G$ coincide, i.e. $F^*(G)=F(G)$.

Write $\alpha$ for a point of $\spaceP$. For a collineation group $H$ of $\spaceP$, write $H_\alpha$ for the stabilizer of $\alpha$.

\section{Nilpotent collineation groups}\label{section:fittingregular}

In order to prove Theorem A we need a well known result of Camina and Praeger. We state a weaker version which is sufficient for our purposes:

\begin{theorem}\label{theorem:campraeg}\cite[Theorem 1]{campraeg2}
Let $G$ act transitively on the points of a projective plane $\spaceP$. Let $N$ be a normal subgroup of $G$. Then $N$ acts faithfully on each of its point orbits.
\end{theorem}

Note that, in particular, Theorem \ref{theorem:campraeg} implies that a minimal normal subgroup of $G$ will act semi-regularly on the points of $\spaceP$. We are now in a position to prove Theorem A.

\begin{proof}
Let $N$ be a Sylow $p$-group of $F$ for some prime $p$ dividing the order of $F$. Suppose that $N$ does not act semi-regularly. If $|N:N_\alpha|=p$, a prime, then, since $N$ acts faithfully, $N<S_p$, the Symmetric group on $p$ letters. But then $|N|=p$ and $N$ is semi-regular which is a contradiction. Thus $|N:N_\alpha|\geq p^2$ and $x>25$.

Observe that the average number of fixed points for non-identity elements of $N$ is
$$(x^2+x+1)\frac{|N_\alpha|-1}{|N|-1}>\frac{x^2+x+1}{2|N:N_\alpha|}.$$

Now $|N:N_\alpha|$ divides into $u^2\pm u+1$. If $|N:N_\alpha|<u^2\pm u+1$ then a non-identity element of $N$ fixes, on average, more than $\frac{3}{2}(u^2-u+1)>u^2+u+1$ fixed points which contradicts our hypothesis.

If $|N:N_\alpha|=u^2\pm u+1$ then, by \cite[p. 11]{ljundggren},
$|N:N_\alpha|=7^3$. 

If $|N:N_\alpha|=7^3=u^2-u+1$ and $|N_\alpha|>7$ then the average number of fixed points for non-identity elements of $N$ is
$$(x^2+x+1)\frac{|N_\alpha|-1}{|N|-1} = 343\times381\frac{|N_\alpha|-1}{7^3|N_\alpha|-1}>\frac{343\times381}{351}>370.$$
Now $x=361$ and an element of order $7$ must fix a multiple of $7$ points. This implies that the most number of points such an element may fix is $357$ which is a contradiction.

If $|N:N_\alpha|=7^3=u^2-u+1$ and $|N_\alpha|=7$ then $|Aut N|$ is not divisible by 127. Thus $F(G)$ must contain a non-trivial Sylow 127-group and $G$ has a normal semi-regular subgroup $P$ of order 127. Furthermore $P$ centralizes $N_\alpha$ and so 
$$|Fix(N_\alpha)|\geq 7\times 127 >x+\sqrt{x}+1.$$ 
This is a contradiction.

If $|N:N_\alpha|=7^3=u^2+u+1$ then $x=324.$ In fact $N=F=F(G)$ since otherwise $F(G)>N$ and a point semi-regular group of order 307 must centralize $N$ which is impossible. 

By Theorem \ref{theorem:campraeg} the centre of $N$ acts semi-regularly on the points of $\spaceP$. This implies that $N_\alpha$ fixes a subplane. Suppose that $|Z(N)|>7$; then $N_\alpha$ must fix a subplane of order $18$. Thus, for $g\in N_\alpha$, $Fix(g)=Fix(N_\alpha)$. Then take $\beta$ a point of $\spaceP$ not in $Fix(N_\alpha)$. We must have $N_\alpha\cap N_\beta=\{1\}$. Thus $|N_\beta|\leq 7^3$ and so $|N|\leq 7^6$. But in this case $307$ does not divide into $|Aut N|$ which is a contradiction.

Suppose, alternatively, that $|Z(N)|=7$. Let $X=N_N(N_\alpha)$. If $X>N_\alpha Z(N)$ then there are at most $7$ $N$-conjugates of $N_\alpha$ and so $N_\alpha$ must fix $0$ or a multiple of 49 points in each orbit of $N$. Hence $N_\alpha$ fixes a Baer subplane and our previous argument can be applied. 

Thus we assume that $X=N_\alpha Z(N)$. Since $N$ is nilpotent $N_N(X)>X$ and so we choose $\alpha$ and $\beta$ points in the same $N$-orbit such that $N_\alpha$ and $N_\beta$ are distinct subgroups of $N_\alpha Z(N)$. Then $|N_\alpha:(N_\alpha\cap N_\beta)|=7$. Furthermore $Y:=N_N(N_\alpha\cap N_\beta)\geq \langle X,n\rangle >X,$ where $n\in N$ such that $\alpha n = \beta,$ hence $Y$ acts on the fixed set of $N_\alpha\cap N_\beta$ with orbits of size a multiple of $49$. We conclude that $N_\alpha\cap N_\beta$ fixes a Baer subplane. 

If $N_\alpha\cap N_\beta$ fixes an entire $N$-orbit then, since $N$ acts faithfully on its point-orbits, $N_\alpha\cap N_\beta=\{1\}$. Thus $|N|=7^4$. Alternatively $N_\alpha\cap N_\beta$ fixes exactly 0 or 49 points in any $N$-orbit. In which case we can find two other points $\gamma$ and $\delta$ in the same $N$-orbit as $\alpha$ such that $N_\gamma\cap N_\delta$ fix a Baer subplane. Then $(N_\alpha\cap N_\beta)\cap(N_\gamma\cap N_\delta)=\{1\}$ and so $|N|\leq 7^8$. Thus in all cases $|N|\leq q^8$ and $307$ does not divide into $|Aut N|$ which is a contradiction.
\end{proof}

\section{Abelian Singer groups}

Throughout this section $S\leq G$ is an abelian Singer group of $\spaceP$.

We record some results of Ho \cite{ho}:

\begin{theorem}\cite[theorem 2]{ho}
An abelian Singer group contained in a soluble collineation group of a finite projective plane is always normal.
\end{theorem}

\begin{lemma}\cite[lemma 4.3]{ho}\label{lemma:centralize}
Let $Q<G$ be a collineation group normalized by $S$. Then $S$ centralizes $Q$ if one of the following holds:
\begin{enumerate}
\item $Q$ is abelian;

\item $|Q|$ is prime to $|S|$;

\item $x=u^2$ and $Q$ is nilpotent.
\end{enumerate}
\end{lemma}

We are now able to give an alternative proof to Theorem C.

\begin{proof}
Suppose that $S$ and $T$ are abelian Singer groups lying in $G$. We need to prove that $S=T$.

If $G$ is soluble then both $S$ and $T$ are normal in $G$ and so lie in $F(G)$. By Lemma \ref{lemma:centralize}, $S$ and $T$ centralize each other. Thus $\langle S,T\rangle$ is abelian and transitive on the points of $\spaceP$. By Theorem \ref{theorem:campraeg}, $|\langle S,T\rangle|\leq v$ and so $S=T$.

If $G$ is not soluble then we may assume that $G$ contains a Baer involution. Hence $x=u^2$ and, by Lemma \ref{lemma:centralize}, $S$ and $T$ both centralize $F(G)$. Since $F(G)=F^*(G)$ this means that $F(G)$ contains both $S$ and $T.$ Now $F(G)$ is soluble and so we can apply the same argument as when $G$ was soluble and conclude that $S=T$.
\end{proof}

\section{Nilpotent Singer groups}\label{section:nilpotentsinger}

Throughout this section $\spaceP$ is a projective plane of order $x=u^2$ and $S\leq G$ is a nilpotent Singer group of $\spaceP$.

\begin{lemma}
$S$ contains $F(G)$.
\end{lemma}
\begin{proof}
Suppose the result does not hold and $S\not\geq F(G)$. Let $P\in Syl_pF(G)$ with $P\not\leq S$. Let $H:=PS=SP$.

Let $Q\in Syl_pH$ such that $Q=PP_1$ where $P_1\in Syl_p S$. Then consider $h^{-1}Qh$ for $h\in H$. We can write $h=sp$ where $p\in P$ and $s\in S$. Then
$$h^{-1}Qh = p^{-1}s^{-1}PP_1sp=P p^{-1}s^{-1}P_1sp = PP_1^p = Q.$$

Thus $Q$ is normal in $H$. But this is a nilpotent group normal in a transitive group hence, by Theorem A, $Q$ is semi-regular. But $|Q|$ does not divide into $x^2+x+1$ and we have a contradiction.
\end{proof}

\begin{corollary}
Any prime dividing into $v$ divides into $|F(G)|$.
\end{corollary}
\begin{proof}
Since $C_G(F(G))\leq F(G)$ we know that if $p$ divides into $|S|$ then $p$ divides into $|F(G)|$. So if $p$ divides into $v$ then $p$ divides into $|F(G)|$.
\end{proof}

We are now in a position to prove Theorem C:

\begin{proof}
Assume, for the sake of contradiction, that $\spaceP$ admits two distinct nilpotent Singer groups $S$ and $T$ and set $G:= Aut \spaceP$. Let $v=p_1^{a_1}\dots p_r^{a_r}$. Then $F(G)=N_1\times\dots\times N_r$, $1\neq N_i\in Syl_{p_i}F(G)$. In a similar way write,

$$S=P_1\times\dots\times P_r, \ \ T=Q_1\times\dots\times Q_r.$$

We will assume, without loss of generality, that $P_1\neq Q_1$. Now $G\leq (N_1\times \dots \times N_r).(A_1\times \dots\times A_r)$ where $A_i$ is a subgroup of the outer automorphism group of $N_i$. Then both $P_i$ and $Q_i$ lie in $(1\times\dots\times 1\times N_i \times 1\times\dots\times 1).(1\times\dots\times 1\times A_i\times 1\times\dots\times 1)$. This implies that 
$$C_G(\langle P_1,Q_1\rangle) \geq P_2\times \dots\times P_r.$$

Now $\langle P_1,Q_1\rangle$ must contain an element $1\neq g$ which fixes a point (since $P_1Q_1$ must have an orbit of size strictly less than $|P_1Q_1|$.) Then we know that
$$\frac{|v|}{|v|_{p_1}} \ \big| \ |Fix \ g|.$$

In fact, consider $\langle Q_1, P_1\rangle$ acting on $\spaceP/ (P_2\times\dots\times P_r)$. Either this is a Frobenius action or there exists $g\in \langle Q_1,P_1\rangle$ that fixes more than one element. But a Frobenius action has a normal Frobenius kernel which must be $Q_1$ and $P_1$. This is a contradiction. Thus

$$\frac{|v|}{|v|_{p_1}} < |Fix \ g|.$$

Now, since $x=u^2$, $v=(x+\sqrt{x}+1)(x-\sqrt{x}+1)$. The bracketed terms are coprime and so $|Fix \ g|\geq 2(x-\sqrt{x}+1)>x+\sqrt{x}+1$. This gives a contradiction.
\end{proof}

\bibliographystyle{amsalpha}
\bibliography{fitting}

\providecommand{\bysame}{\leavevmode\hbox to3em{\hrulefill}\thinspace}
\providecommand{\MR}{\relax\ifhmode\unskip\space\fi MR }
\providecommand{\MRhref}[2]{%
  \href{http://www.ams.org/mathscinet-getitem?mr=#1}{#2}
}
\providecommand{\href}[2]{#2}
\begin{thebibliography}{Dem97}

\bibitem[CP93]{campraeg2}
Alan~R. Camina and Cheryl~E. Praeger, \emph{Line-transitive automorphism groups
  of linear spaces}, Bull. London Math. Soc. \textbf{25} (1993), 309--315.

\bibitem[Dem97]{dembov}
P.~Dembowski, \emph{Finite geometries}, Springer-Verlag, 1997.

\bibitem[Gil]{gill}
Nick Gill, \emph{Transitive projective planes}, Submitted.

\bibitem[Ho98]{ho}
C.~Y. Ho, \emph{Finite projective planes with abelian transitive collineation
  groups}, J. Algebra \textbf{208} (1998), 533--550.

\bibitem[Kan87]{kantor}
W.~Kantor, \emph{Primitive permutation groups of odd degree, and an application
  to finite projective planes}, J. Algebra \textbf{106} (1987), 15--45.

\bibitem[Lju43]{ljundggren}
W.~Ljunggren, \emph{Einige bemerkungen {\"u}ber die {D}arstellung ganzer
  {Z}ahlen durch bin{\"a}re kubische {F}ormen mit positiver {D}iskriminante},
  Acta. Math. \textbf{74} (1943), 1--21.

\bibitem[Wag59]{wagner}
A.~Wagner, \emph{On perspectivities of finite projective planes}, Math. Z.
  \textbf{71} (1959), 113--123.

\end{thebibliography}

\end{document}